\documentclass[12pt,a4paper]{amsart}

\usepackage{amsmath}
\usepackage{amsthm}
\usepackage{amssymb}
\usepackage{amscd}
\usepackage[all]{xy}

\title[The Nef Curve Cone Theorem Revisited]
{The Nef Curve Cone Theorem Revisited}
\author{Qihong Xie}
\date{2005/01/20}
\thanks{}
\address{Graduate School of Mathematical Sciences, University of Tokyo, 
3-8-1 Komaba, Meguro, Tokyo 153-8914, Japan}
\email{xqh@ms.u-tokyo.ac.jp}

\theoremstyle{plain}
\newtheorem{prop}{Proposition}[section]
\newtheorem{lem}[prop]{Lemma}
\newtheorem{thm}[prop]{Theorem}
\newtheorem{cor}[prop]{Corollary}
\newtheorem{conj}[prop]{Conjecture}

\theoremstyle{definition}
\newtheorem{defn}[prop]{Definition}
\newtheorem*{ack}{Acknowledgements}

\theoremstyle{remark}
\newtheorem{rem}[prop]{Remark}

\newcommand{\Q}{\mathbb Q}
\newcommand{\R}{\mathbb R}
\newcommand{\C}{\mathbb C}
\newcommand{\Z}{\mathbb Z}
\newcommand{\N}{\mathbb N}

\newcommand{\PP}{\mathbb P}

\newcommand{\EE}{\mathcal E}

\newcommand{\LL}{\mathcal L}

\newcommand{\TT}{\mathcal T}

\newcommand{\NE}{\overline{NE}}
\newcommand{\NM}{\overline{NM}}
\newcommand{\Eff}{\overline{Eff}}
\newcommand{\ep}{\varepsilon}
\newcommand{\dra}{\dashrightarrow}
\newcommand{\ra}{\rightarrow}
\newcommand{\inter}[2]{\langle #1,#2 \rangle}
\newcommand{\Int}{\mathop{\rm Int}\nolimits}
\newcommand{\codim}{\mathop{\rm codim}\nolimits}

\begin{document}

\begin{abstract}
We revisit the nef curve cone theorem 
and use it to reprove the pseudo-effectivity of the second 
Chern classes for terminal weak $\Q$-Fano varieties.
\end{abstract}

\maketitle

\setcounter{section}{0}
\section{Introduction}\label{S1}

The purpose of this note is to revisit the nef curve cone theorem 
due to \cite{ba}. We state this theorem in the higher dimensional 
log case and give a proof almost identical to the original one. 
As an important application, we use this theorem to reprove 
the pseudo-effectivity of the second Chern classes for 
terminal weak $\Q$-Fano varieties, which was, in fact, proved 
in \cite{kmmt}.

There are three reasons to redo this interesting work. 
First, we would like to give the formal statement and proof of 
the nef curve cone theorem in the higher dimensional log case. 
Second, this theorem in the three-dimensional log case was 
claimed in \cite{kmmt}, however some necessary explanations 
were absent. Third, some inaccurate points and typographical 
errors in \cite{ba} and \cite{kmmt} are clarified in this note.

Before stating the main results, we need the following important 
definitions of cones in birational geometry (cf.\ \cite{la}).

Let $X$ be a normal proper variety. Let $N^1(X)$ (resp.\ $N_1(X)$) 
be the $\R$-vector space generated by all Cartier divisors (resp.\ 
all irreducible proper curves) on $X$ modulo numerical equivalence. 
The intersection number $N_1(X)\times N^1(X)\ra \R$, 
$(C,D)\mapsto \inter{C}{D}=C\cdot D$ defines a nondegenerate pairing 
between $N_1(X)$ and $N^1(X)$.

The effective divisor cone $Eff(X)\subset N^1(X)$ 
is defined to be a cone generated by all effective Cartier divisors 
on $X$. Its closure $\Eff(X)$ is called the pseudo-effective 
divisor cone. A $\Q$-Cartier divisor $D$ is pseudo-effective if and 
only if the Iitaka D-dimension $\kappa(X,D)\geq 0$. 
Let $\NM(X)\subset N_1(X)$ be the dual cone of $\Eff(X)$, 
which is called the nef curve cone.

$NE(X)\subset N_1(X)$ is defined to be a cone generated by all 
irreducible proper curves on $X$. Its closure $\NE(X)$ is 
called the Kleiman-Mori cone, which plays an essential role in the 
Minimal Model Theory. Let $B$ be an effective $\Q$-divisor on $X$ 
such that $K_X+B$ is $\Q$-Cartier. 
Let $H$ be an ample Cartier divisor on $X$, 
$\ep$ a sufficiently small positive number. 
Denote $\NE_\ep(X)=\{ z\in\NE(X)$ $|$ $\inter{z}{K_X+B+\ep H}\geq 0 \}$. 
Let $Nef(X)\subset N^1(X)$ be the dual cone of 
$\NE(X)$, which is called the nef divisor cone and indeed 
generated by all nef Cartier divisors on $X$.

The big divisor cone $Big(X)\subset N^1(X)$ is defined to be a 
cone generated by all big Cartier divisors on $X$. 
If $X$ is projective, we define the ample divisor cone 
$Amp(X)\subset N^1(X)$ to be a cone generated by all ample 
Cartier divisors on $X$. It is easy to see that both $Big(X)$ and 
$Amp(X)$ are open cones, and the following relations hold:
\[ Big(X)=\Int(\Eff(X)),\quad \overline{Amp}(X)=Nef(X). \]

The following is the nef curve cone theorem.

\begin{thm}\label{1.1}
Let $X$ be a normal projective $\Q$-factorial variety of dimension $n$, 
$B$ an effective $\Q$-divisor on $X$ such that $(X,B)$ is 
Kawamata log terminal (klt, for short). 
Assume that $K_X+B$ is not pseudo-effective. 
Let $H$ be an ample Cartier divisor on $X$, 
$\ep$ a sufficiently small positive number. 
Then up to the Log Minimal Model Conjecture and 
Conjecture \ref{1.3} on dimension $n$, 
there are a positive integer $r$ depending only on $\ep$ and 
integer curves $l_i\in\NM(X)\setminus\NE_\ep(X)$ for 
$1\leq i\leq r$, such that
\[ \NM(X)+\NE_\ep(X)=\sum_{i=1}^r\R_+[l_i]+\NE_\ep(X). \]
\end{thm}

\begin{conj}[Log Minimal Model Conjecture]\label{1.2}
Let $X$ be a normal projective $\Q$-factorial variety of dimension $n$, 
$B$ an effective $\Q$-divisor on $X$ such that $(X,B)$ is klt. 
Then there exists a sequence of birational maps
\[(X=X_0,B=B_0)\dra (X_1,B_1)\dra 
\cdots\dra (X_k,B_k)\]
where $B_i$ is the strict transform of $B$ on $X_i$ 
and $\dra$ is either an extremal divisorial contraction or 
a log flip, such that $(X_k,B_k)$ is either a log minimal model or a 
log Mori fiber space.
\end{conj}

\begin{conj}\label{1.3}
Notation and assumptions as above. Furthermore assume that $K_X+B$ is nef. 
Then the Iitaka D-dimension $\kappa(X,K_X+B)\geq 0$.
\end{conj}

For the Log Minimal Model Conjecture, it is sufficient to prove 
the existence of log flips and the termination of log flips 
(cf.\ \cite{kmm}). 
The existence of log flips is known on dimension at most four 
(cf.\ \cite{sh92,sh03}), and the termination of log flips is known 
on dimension at most three (cf.\ \cite{ka92}).

Instead of Conjecture \ref{1.3}, in fact, we may put forward 
a stronger one called the Log Abundance Conjecture, 
which was proved on dimension at most three (cf.\ \cite{kmm94}).

\begin{rem}\label{1.4}
Theorem \ref{1.1} is also valid for divisorial log terminal pairs 
provided that the corresponding two conjectures are true for dlt pairs.
\end{rem}

As an application of Theorem \ref{1.1}, we reprove the following

\begin{thm}\label{1.5}
Let $X$ be a terminal weak $\Q$-Fano variety of dimension $n$, i.e.\ 
$X$ has only terminal singularities and $-K_X$ is nef and big. 
Then up to the Log Minimal Model Conjecture and 
Conjecture \ref{1.3} on dimension $n$, for any ample 
Cartier divisors $H_1,\cdots,H_{n-2}$ on $X$, we have 
\[ c_2(X)\cdot H_1\cdots H_{n-2}\geq 0. \]
\end{thm}

In \S \ref{S2} and \S \ref{S3}, we will prove Theorems \ref{1.1} and 
\ref{1.5} respectively. 
In what follows, we assume that the Log Minimal Model Conjecture 
and Conjecture \ref{1.3} hold on dimension $n$. We will use freely 
the results on the Minimal Model Theory, and refer the reader to 
\cite{kmm, km} for more details.

We work over the complex number field $\C$.

\begin{ack}
I would like to express my gratitude to Professor Hiromichi Takagi 
for his useful suggestions and comments.
\end{ack}

\section{Proof of Theorem \ref{1.1}}\label{S2}

\begin{defn}\label{2.1}
Let $\varphi: X\dra Z$ be a birational map between normal proper 
varieties. $\varphi$ is called an ICO-map if the birational morphism 
$\varphi^{-1}: Z\setminus Z_0\ra\varphi^{-1}(Z\setminus Z_0)$ is 
isomorphic in codimension one, where $Z_0$ is the indeterminancy locus 
of $\varphi^{-1}$ with codimension at least two.
\end{defn}

For such an ICO-map $\varphi$, we have one natural injective map 
$\varphi^{1*}: N^1(Z)\ra N^1(X)$, and by duality, the other 
natural injective map $\varphi_1^{*}: N_1(Z)\ra N_1(X)$. 
Note that $\varphi_1^{*}$ is not the strict transform of 1-cycles, 
since $\varphi_1^{*}(N_1(Z))\subset N_1(X)$ is always orthogonal to 
the kernel of the direct image map $N^1(X)\ra N^1(Z)$.
It is easy to see that the composition of any two ICO-maps is also an 
ICO-map, and that extremal divisorial contractions and log flips are 
ICO-maps.

In \S\ref{S2}, we fix the notation and assumptions as in Theorem \ref{1.1}.

\begin{defn}\label{2.2}
Let $H$ be an ample Cartier divisor on $X$. 
The following three numbers are improtant for technical reasons.

$\alpha_X(H):=\inf\{p/q\in\Q$ $|$ $p\in\Z$, $q\in\N$, $pH+q(K_X+B)\in 
Big(X) \}$;

$\sigma_X(H):=\sup\{t\in\R$ $|$ $H+t(K_X+B)\in \Eff(X) \}$;

$\tau_X(H):=\sup\{t\in\R$ $|$ $H+t(K_X+B)\in Nef(X) \}$.
\end{defn}

\begin{lem}\label{2.3}
$\sigma_X(H)=\alpha_X(H)^{-1}$. If for $m\gg 0$, 
$m(H+\tau_X(H)(K_X+B))$ defines a morphism from 
$X$ to $Z$ with $\dim Z<\dim X$, then $\sigma_X(H)=\tau_X(H)$.
\end{lem}

\begin{proof}
Since every big divisor is pseudo-effective, we have $\sigma_X(H)\geq
\alpha_X(H)^{-1}$. The $\Q$-divisor $(1+\ep_1)H+(\sigma_X(H)-\ep_2)
(K_X+B)$ is big for suitable positive number $\ep_1,\ep_2$. 
Hence $\alpha_X(H)^{-1}\geq(\sigma_X(H)-\ep_2)/(1+\ep_1)$. 
Let $\ep_1,\ep_2\to 0$, we obtain $\sigma_X(H)=\alpha_X(H)^{-1}$.

Moreover, If for $m\gg 0$, $m(H+\tau_X(H)(K_X+B))$ 
defines a morphism from $X$ to $Z$ with $\dim Z<\dim X$, 
then $H+\tau_X(H)(K_X+B)$ is not big, hence $\tau_X(H)\geq
\alpha_X(H)^{-1}$. Since every nef divisor is pseudo-effective, 
we always have $\tau_X(H)\leq\sigma_X(H)=\alpha_X(H)^{-1}$.
\end{proof}

\begin{thm}\label{2.4}
Assume that $K_X+B$ is not pseudo-effective. Then for any nonempty open 
subset $U\subset Amp(X)$, there are a variety $S$ of dimension $n-1$ 
with a dominant rational map $\psi: \PP^1\times S\dra X$, 
an ICO-map $\varphi: X\dra Z$, and a nonempty open subset 
$V\subset U$ such that

(i) for every general point $s\in S$, 
$0<\inter{\psi_*(\PP^1\times s)}{-(K_X+B)}\leq 2n$;

(ii) $\sigma_X(L)$ is a rational number for any $\Q$-Cartier divisor 
$L\in V$;

(iii) for $m\gg 0$ the linear system $|m(L+\sigma_X(L)(K_X+B))|$ is 
nonempty for any $\Q$-Cartier divisor $L\in V$, and 
$\inter{\varphi^*_1\varphi_*\psi_*(\PP^1\times s)}{L+\sigma_X(L)(K_X+B)}=0$ 
for the generic point $s\in S$.
\end{thm}

\begin{proof}
Since $U$ is an open subset consisting of ample $\R$-Cartier divisors, 
we can find an ample $\Q$-Cartier divisor $L_0\in U$ such that the 
divisor $D=L_0+\tau_0(K_X+B)$ defines an extremal ray $R=\NE(X)\cap D^
\bot$, where $\tau_0=\tau_X(L_0)$.

For any $\Q$-Cartier divisor $L\in U$, denote $D(L)=L+\tau_X(L)(K_X+B)$. 
Without loss of generality, we can assume that $\NE(X)\cap D^\bot=R$ 
for all $\Q$-Cartier divisor $L\in U$, 
otherwise we can replace $U$ by some smaller neighbourhood of $L_0$. 
Let $f: X\ra Y$ be the corresponding extremal contraction of $R$.

Up to the Log Minimal Model Conjecture and Conjecture \ref{1.3}
on dimension $n$, since $K_X+B$ is not pseudo-effective, after taking 
a finite number of divisorial contractions or log flips, at last we reach 
to a log Mori fiber space. The main idea of the proof is to verify 
Theorem \ref{2.4} for the log Mori fiber space $f: X\ra Y$, and 
use induction to deal with the divisorial contraction case and the 
log flip case.

{\it Case I: $f$ is a log Mori fiber space.}

By Lemma \ref{2.3} and the rationality theorem, we have $\sigma_X(L)=
\tau_X(L)$ is a rational number for all $\Q$-Cartier divisor 
$L\in U$. On the other hand, $X$ is covered by a family of rational 
curves whose general member can be chosen to be an 
extremal rational ray $l$ with length $\inter{l}{-(K_X+B)}\leq 2n$ 
(cf.\ \cite[Theorem 1]{ka91}), and for $m\gg 0$, the linear system 
$|mD(L)|$ is base point free and vertical to $l$. Thus we can take a 
$(n-1)$-dimensional family $\psi: \PP^1\times S\dra X$ with 
$\psi_*(\PP^1\times s)=l$, $V=U$, $Z=X$ and $\varphi=id$ to complete 
the proof of Case I.

{\it Case II: $f$ is a divisorial contraction.}

$(Y,B_Y=f_*B)$ is also $\Q$-factorial and klt. We may write 
\[ K_X+B=f^*(K_Y+B_Y)+aE, \]
where $a$ is a positive rational number and $E$ is the exceptional 
divisor of $f$. Let $F: U\ra N^1(Y)$ be the affine linear map 
defined by 
\[ F(W)=f_*(W)+\tau_0(K_Y+B_Y). \]
Since $F$ is an open map, we get a subset $F(U)\subset N^1(Y)$ containing 
the ample $\Q$-Cartier divisor $F(L_0)=f_*(L_0)+\tau_0(K_Y+B_Y)=
f_*(L_0+\tau_0(K_X+B))=f_*(D(L_0))$. Let $U_0=F(U)\cap Amp(Y)$. 
Then $F(L_0)\in U_0\neq\emptyset$.

{\it Claim 1. $\sigma_Y(F(L))=\sigma_X(L)-\tau_0$ for any $\Q$-Cartier 
divisor $L\in V:=F^{-1}(U_0)\cap U$.}

{\it Proof of Claim 1.} Since $f_*\Eff(X)=\Eff(Y)$, we have
\begin{eqnarray}
\sigma_Y(F(L)) & = & \sigma_Y(f_*(L+\tau_0(K_X+B))) \nonumber \\
& \geq & \sigma_X(L+\tau_0(K_X+B)) \nonumber \\
& = & \sigma_X(L)-\tau_0. \nonumber
\end{eqnarray}
The divisor $G(L)=f^*(F(L)+\sigma_Y(F(L))(K_Y+B_Y))$ is the total 
transform of a pseudo-effective divisor on $Y$, so $G(L)$ is also 
pseudo-effective. On the other hand, we may write 
\[ G(L)=L+(\tau_0+\sigma_Y(F(L)))(K_X+B)+bE, \]
where $b$ is a rational number. Let $l\in R$ be a curve, then 
$\inter{l}{G(L)}=0$. Since $\inter{l}{(K_X+B)}<0$, $\inter{l}{E}<0$, 
$\tau_0+\sigma_Y(F(L))\geq\sigma_X(L)\geq\tau_X(L)$, and 
$\inter{l}{L+\tau_X(L)(K_X+B)}=0$, we have
\[ -b\inter{l}{E}=\inter{l}{L+(\tau_0+\sigma_Y(F(L)))(K_X+B)}\leq 0, \]
so $b\leq 0$. Hence $G(L)-bE=L+(\tau_0+\sigma_Y(F(L)))(K_X+B)$ is 
pseudo-effective and $\sigma_X(L)\geq\tau_0+\sigma_Y(F(L))$, which 
completes the proof of Claim 1.

By the induction hypothesis for $Y$, we can find a $(n-1)$-fold $S$ 
with a dominant rational map $\psi_0: \PP^1\times S\dra Y$, 
an ICO-map $\chi: Y\dra Z$ and an open subset $V_0\subset U_0$, 
such that the conditions (i)-(iii) are satisfied.

Let $\psi=f^{-1}\circ\psi_0: \PP^1\times S\dra X$ be the strict 
transform of the family $\psi_0$ via $f$. Since all curves of the 
family $\psi$ do not lie in $E$, hence $\inter{\psi_*(\PP^1\times s)}
{E}\geq 0$ and 
\begin{eqnarray}
\inter{\psi_*(\PP^1\times s)}{-(K_X+B)} & \leq & 
\inter{\psi_*(\PP^1\times s)}{-f^*(K_Y+B_Y)} \nonumber \\
& = &\inter{\psi_{0*}(\PP^1\times s)}{-(K_Y+B_Y)}\leq 2n. \nonumber
\end{eqnarray}

It follows from Claim 1 that $\sigma_X(L)$ is a rational number for 
any $\Q$-Cartier divisor $L\in V$.

For $m\gg 0$, the direct image of the linear system $|m(L+\sigma_X(L)
(K_X+B))|$ is the nonempty linear system $|m(F(L)+\sigma_Y(F(L))
(K_Y+B_Y))|$. Let $\varphi=\chi\circ f: X\dra Z$. 
Since $\inter{\varphi^*_1\varphi_*\psi_*(\PP^1\times s)}{E}=0$ and 
$\inter{\chi^*_1\chi_*\psi_{0*}(\PP^1\times s)}
{F(L)+\sigma_Y(F(L))(K_Y+B_Y)}=0$, it is easy to see that 
$\inter{\varphi^*_1\varphi_*\psi_*(\PP^1\times s)}{L+\sigma_X(L)
(K_X+B)}=0$, which completes the proof of Case II.

{\it Case III: $f$ is a flipping contraction.}

Up to the Log Minimal Model Conjecture, there exists a log flip 
$tr_R: X\dra X^+$ such that 
\[ \xymatrix{
X \ar@{.>}[rr]^{tr_R} \ar[dr]_f &   & X^+ \ar[dl]^{f^+} \\
 & Y & } \]

(1) $tr_R$ is isomorphic in codimension one;

(2) $(X^+,B^+=tr_{R*}(B))$ is $\Q$-factorial klt, and $K_{X^+}+B^+$ is 
$f^+$-ample.

Thus $tr_R$ yields an isomorphism $\beta_R: \Eff(X)\ra\Eff(X^+)$. 
Since $f_*(L_0+\tau_0(K_X+B))$ is an ample $\Q$-Cartier divisor on $Y$, 
we have that for $m\gg 0$, $\beta_R(L_0+(\tau_0+1/m)(K_X+B))$=
$f^{+*}f_*(L_0+\tau_0(K_X+B))+1/m(K_{X^+}+B^+)$ is an ample 
$\Q$-Cartier divisor on $X^+$.

Let $F: U\ra N^1(X^+)$ be the affine linear map defined by 
\[ F(W)=tr_{R*}(W+(\tau_0+1/m)(K_X+B)). \]
Since $F$ is an open map, we obtain a subset $F(U)\subset N^1(X^+)$ 
containing the ample $\Q$-Cartier divisor $F(L_0)$. Let $U_0=F(U)\cap 
Amp(X^+)$. Then $F(L_0)\in U_0\neq\emptyset$.

By the induction hypothesis for $X^+$, we can find a $(n-1)$-fold $S$ 
with a dominant rational map $\psi_0: \PP^1\times S\dra X^+$, 
an ICO-map $\chi: X^+\dra Z$ and an open subset $V_0\subset U_0$, 
such that the conditions (i)-(iii) are satisfied. 
Let $\psi: \PP^1\times S\dra X$ be the strict transform of the family 
$\psi_0$ via $tr_R$, and let $V=F^{-1}(V_0)\cap U$.

{\it Claim 2. $\inter{\psi_*(\PP^1\times s)}{-(K_X+B)}\leq 
\inter{\psi_{0*}(\PP^1\times s)}{-(K_{X^+}+B^+)}$ for the generic 
point $s\in S$}.

{\it Proof of Claim 2.} Let $\widetilde{X}$ be a common resolution of 
$X$ and $X^+$.
\[ \xymatrix{
  & \ar[dl]_g \widetilde{X} \ar[dr]^{g^+} & \\
X &  & X^+} \]
We may write
\begin{eqnarray}
K_{\widetilde{X}} & = & g^*(K_X+B)+\sum a_iE_i, \nonumber \\
K_{\widetilde{X}} & = & g^{+*}(K_{X^+}+B^+)+\sum b_iE_i. \nonumber
\end{eqnarray}
Then we have $a_i\leq b_i$ for all $i$. Let $\widetilde{\psi}$ be the 
strict transform on $\widetilde{X}$ of the family $\psi$ via $g$. Then 
$\inter{\widetilde{\psi}_*(\PP^1\times s)}{E_i}\geq 0$ for all $i$ since 
the family $\widetilde{\psi}$ covers a dense subset of $\widetilde{X}$. 
Thus
\begin{eqnarray}
& & \inter{\psi_*(\PP^1\times s)}{-(K_X+B)} \nonumber \\
& = & \inter{\widetilde{\psi}_*(\PP^1\times s)}{-K_{\widetilde{X}}}+
\sum a_i\inter{\widetilde{\psi}_*(\PP^1\times s)}{E_i} \nonumber \\
& \leq & \inter{\widetilde{\psi}_*(\PP^1\times s)}{-K_{\widetilde{X}}}+
\sum b_i\inter{\widetilde{\psi}_*(\PP^1\times s)}{E_i} \nonumber \\
& = & \inter{\psi_{0*}(\PP^1\times s)}{-(K_{X^+}+B^+)}. \nonumber
\end{eqnarray}
This completes the proof of Claim 2, which guarantees the condition (i) 
for $(X,B)$.

Since the isomorphism $\beta_R$ gives $\sigma_{X^+}(F(L))=\sigma_X(L)-\tau_0
-1/m$, we have the rationality of $\sigma_X(L)$ for any $\Q$-Cartier 
divisor $L\in V$.

For $m\gg 0$, the direct image of the linear system $|m(L+\sigma_X(L)
(K_X+B))|$ is the nonempty linear system $|m(F(L)+\sigma_{X^+}(F(L))
(K_{X^+}+B^+))|$. Let $\varphi=\chi\circ tr_R: X\dra Z$. Then (iii) is 
satisfied for $(X,B)$, which completes the proof of Case III.
\end{proof}

\begin{cor}\label{2.5}
Assume that $K_X+B$ is not pseudo-effective. 
Let $H$ be an ample $\R$-Cartier divisor on $X$. Then there exist a 
$(n-1)$-fold $S$ with a dominant rational map $\psi: \PP^1\times S\dra X$, 
and an ICO-map $\varphi: X\dra Z$ such that $\sigma_X(H)=\inter{C}{H}/
\inter{C}{-(K_X+B)}$, where $C=\varphi^*_1\varphi_*\psi_*(\PP^1\times s)$, 
$0<\inter{\psi_*(\PP^1\times s)}{-(K_X+B)}\leq\inter{C}{-(K_X+B)}\leq 2n$ 
for the generic point $s\in S$.
\end{cor}

\begin{proof}
By Theorem \ref{2.4}, we can find a sequence of ample $\Q$-Cartier 
divisors $\{ H_i \}$ and $(n-1)$-dimensional covering families 
$\{ \psi_i: \PP^1\times S_i\dra X \}$ and ICO-maps $\varphi_i: X\dra Z_i$ 
such that $\lim_{i\to\infty}H_i=H$ and $\sigma_X(H_i)=\inter{C_i}{H_i}/
\inter{C_i}{-(K_X+B)}$, where $C_i=(\varphi_i)^*_1(\varphi_i)_*(\psi_i)_*
(\PP^1\times s)$.
Since $\sigma_X(\cdot)$ is a continuous function on the divisors, 
we have $\lim_{i\to\infty}\sigma_X(H_i)=\sigma_X(H)$, 
hence $\inter{C_i}{H_i}/\inter{C_i}{-(K_X+B)}$ are bounded. 
Note that $\inter{C_i}{-(K_X+B)}$ are positive and 
bounded from above by $2n$. So $\inter{C_i}{H_i}$, hence 
$\inter{C_i}{H}$ are bounded. Note also that the curves $C_i$ 
are the integer points in $\NM(X)$, so there exist only 
finitely many integer points in this bounded subset of $\NM(X)$.

Let $C$ be an element which repeats in the sequence $\{ C_i \}$ 
infinitely many times. Then $\sigma_X(H)=\inter{C}{H}/
\inter{C}{-(K_X+B)}$, and the family $\psi$ in $\{ \psi_i \}$ 
corresponding to $C$ is the desirable one.
\end{proof}

\begin{cor}\label{2.6}
Assume that $K_X+B$ is not pseudo-effective. 
Let $H$ be an ample Cartier divisor on $X$. Then $\sigma_X(H)$ is a 
rational number.
\end{cor}

\begin{proof}
By Corollary \ref{2.5}, there exists an integer curve $C$ such that 
\[ \sigma_X(H)=\inter{C}{H}/\inter{C}{-(K_X+B)}. \]
Note that all the intersection numbers are rational, so we are done.
\end{proof}

\begin{defn}\label{2.7}
A half line $R=\R_+[l]\subset\NM(X)$ is called a coextremal ray if 
the following conditions hold:

(i) $l\not\in\NE_0(X)=\{ z\in\NE(X)$ $|$ $\inter{z}{K_X+B}\geq 0 \}$;

(ii) if $z_1,z_2\in\NM(X)+\NE_0(X)$ and $z_1+z_2\in R$, then $z_1,z_2\in R$.
\end{defn}

\begin{proof}[Proof of Theorem \ref{1.1}]
Let $P(\ep)=\{$ integer curves $l\in\NM(X)$ $|$ $\inter{l}{-(K_X+B)}\leq 2n$, 
$l\not\in\NE_\ep(X)\}$. All integer curves $l\in P(\ep)$ have degree 
$\inter{l}{H}\leq 2n/\ep$, hence form a finite set in $\NM(X)$, which 
is denoted by $\{l_1, \cdots, l_r\}$.

We shall prove that 
\[ T(\ep):=\R_+[l_1]+\cdots+\R_+[l_r]+\NE_\ep(X)=\NM(X)+\NE_\ep(X). \]
Obviously, $T(\ep)\subset\NM(X)+\NE_\ep(X)$. Assume that 
$T(\ep)\neq\NM(X)+\NE_\ep(X)$. Then by the separation theorem for 
convex sets, there exists a $\Q$-Cartier divisor $D\in N^1(X)$ such that 

(1) $\inter{z}{D}\geq 0$ for arbitrary $z\in\NM(X)+\NE_\ep(X)$;

(2) there exists a $0\neq z_0\in\NM(X)$ such that $\inter{z_0}{D}=0$;

(3) for any $z\in T(\ep)\setminus \{0\}$, the number $\inter{z}{D}/
\inter{z}{H}$ is positive and separated from zero by an absolute constant.

Let $G$ be the convex cone of $N^1(X)$ generated by $D$ and $-(K_X+B)$.

{\it Claim 1. $Amp(X)\cap G\neq \{0\}$. }

{\it Proof of Claim 1.} Assume that $Amp(X)\cap G=\{0\}$. Then by the 
separation theorem, there exists an element $u\in N_1(X)$ such that 
$\inter{u}{Amp(X)}\geq 0$ and $\inter{u}{G\setminus \{0\}}<0$. From 
the first inequality we see that $u\in\NE(X)$, and from the second one 
we have $\inter{u}{D}<0$ and $\inter{u}{-(K_X+B)}<0$. Hence 
$u\in\NE_\ep(X)$ but $\inter{z}{D}<0$, a contradiction to (1), which 
completes the proof of Claim 1.

Thus there exists an ample $\Q$-Cartier divisor 
$D_0=\lambda_0D-\mu_0(K_X+B)$ for some $\lambda_0>0,\mu_0>0$, 
since if $\lambda_0=0$, then $-(K_X+B)$ is ample and we can change 
$\lambda_0,\mu_0$ slightly such that $\lambda_0>0,\mu_0>0$. 
By Corollary \ref{2.5}, there exist a $(n-1)$-fold $S$ with a 
dominant rational map $\psi: \PP^1\times S\dra X$, 
and an ICO-map $\varphi: X\dra Z$ such that 
$\sigma_X(D_0)=\inter{C}{D_0}/\inter{C}{-(K_X+B)}$, 
where $C=\varphi^*_1\varphi_*\psi_*(\PP^1\times s)$, 
s is the generic point of $S$ and 
$0<\inter{\psi_*(\PP^1\times s)}{-(K_X+B)}\leq\inter{C}{-(K_X+B)}\leq 2n$.

By definition, we have $C\in T(\ep)$. 
It follows from (1) that $D$ is pseudo-effective, hence 
$\sigma_X(D_0)\geq\mu_0$. For $0<\delta\ll 1$, we have 
\[ D_0+(\mu_0+\delta)(K_X+B)=\frac{\lambda_0(\mu_0+\delta)}{\mu_0}D-
\frac{\delta}{\mu_0}D_0 \]
is not pseudo-effective by (2). Thus $\sigma_X(D_0)=\mu_0$. 
But by some easy computation, we have $\inter{C}{D}=0$, which 
contradicts to (3).
\end{proof}

\begin{rem}\label{2.8}
With notation and assumptions as in Theorem \ref{1.1}. Assume that 
at last we obtain a log Mori fiber space $g: X'\ra Z$ by running 
$(K_X+B)$-MMP. Let $f: X\dra X'$ be the natural birational map and 
$U'=\{x'\in X'$ $|$ $f^{-1}$ is an isomorphism at $x'\}$. Since 
$\codim(X'\setminus U')\geq 2$, we can take a projective curve 
$l'\subset U'$ contained in a general fiber of $g$. Let $l$ be the 
strict transform of $l'$ on $X$. Then $\R_+[l]$ is a coextremal ray for 
some $\ep$. Conversely, every coextremal ray is obtained by this way.
\end{rem}

\begin{cor}\label{2.9}
With notation and assumptions as in Theorem \ref{1.1}. Assume further 
that $-(K_X+B)$ is ample. Then $\NM(X)=\sum_{i=1}^r\R_+[l_i]$ 
is a polyhedral cone.
\end{cor}

\begin{proof}
Take $H=-(K_X+B), \ep=1/2$ and apply Theorem \ref{1.1}.
\end{proof}

\section{Proof of Theorem \ref{1.5}}\label{S3}

In \S\ref{S3}, we fix the notation and assumptions as in Theorem \ref{1.5}.

\begin{defn}\label{3.1}
Let $X$ be a variety with only terminal singularities, $U$ the smooth 
locus of $X$. Then $\codim(X\setminus U)\geq 3$. The second Chern class 
$c_2(X)$ of $X$ is defined as follows:
\[ c_2(X)=c_2(\TT_X):=c_2(\TT_X|_U)\in A^2(U)\cong A^2(X), \]
where $\TT_X$ is the tangent sheaf on $X$.
\end{defn}

In order to prove the pseudo-effectivity of $c_2(X)$, we make 
use of the following criterion (cf.\ \cite[Theorem 6.1]{mi}).

\begin{thm}\label{3.2}
Let $X$ be a normal projective variety of dimension $n$ which is 
smooth in codimension two. Let $H_1,\cdots,H_{n-1}$ be ample Cartier 
divisors on $X$. Let $\EE$ be a torsion free sheaf on $X$ such that

(1) $c_1(\EE)$ is a nef $\Q$-Cartier divisor;

(2) $\EE$ is generically $(H_1,\cdots,H_{n-1})$-semipositive, i.e. 
\[ c_1(\LL)\cdot H_1\cdots H_{n-1}\geq 0, \]

for every torsion free quotient sheaf $\LL$ of $\EE$.

\noindent Then $c_2(\EE)\cdot H_1\cdots H_{n-2}\geq 0$.
\end{thm}

\begin{defn}\label{3.3}
Let $X$ be a variety and $C$ a rational curve contained in the smooth 
locus of $X$. $C$ is called a free rational curve if $\TT_X|_C$ is 
semipositive.
\end{defn}

We recall the following proposition (cf.\ \cite[Corollary 1.3]{kmm92}), 
which says that there exists a free rational curve in any uniruled 
variety.

\begin{prop}\label{3.4}
Let $X$ be a variety covered by a family of rational curves $\{C\}$ 
such that $C$ is contained in the smooth locus of $X$. Then a general 
member $C$ is a free rational curve.
\end{prop}

The existence of free rational curves and the nef curve cone theorem 
are the keypoints in the proof of Theorem \ref{1.5}.

\begin{proof}[Proof of Theorem \ref{1.5}]
By Kodaira's lemma, we can find an effective $\Q$-divisor $B$ on $X$, 
such that $(X,B)$ is terminal and $-(K_X+B)$ is ample. Hence 
$\NM(X)=\sum_{i=1}^r\R_+[l_i]$ by Corollary \ref{2.9}.

Let $\EE=\TT_X$. Then $c_1(\TT_X)=-K_X$ is nef. Hence by Theorem \ref{3.2}, 
it suffices to prove that $\EE$ is generically $(H_1,\cdots,H_{n-1})$
-semipositive. Since $H_1\cdots H_{n-1}\in\NM(X)$, we have only to 
prove that $c_1(\LL)\cdot l_i\geq 0$ for each $1\leq i\leq r$ and each 
surjection $\EE\twoheadrightarrow\LL$ to a torsion free sheaf $\LL$.

We use the description of coextremal rays in Remark \ref{2.8} and fix a 
coextremal ray $l\in\{l_1,\cdots,l_r\}$. Let $U$ be the open subset 
of $X$ corresponding to $U'$. We extend $\TT_X|_U\ra\LL|_U/({\rm tor})$ 
to $\TT_{X'}\ra \LL'$ via $U\cong U'\subset X'$. 
Since $X\dra X'$ is an isomorphism on $U\supset l$, we have 
$c_1(\LL)\cdot l=c_1(\LL')\cdot l'$.

By construction, we may assume that $l'$ is contained in a fiber of 
$X'\ra Z$ and is a general member in a covering family of rational 
curves on $X'$. Hence $l'$ is a free rational curve on $X'$. 
Note that $\TT_{X'}|_{l'}\ra \LL'|_{l'}$ has a finite cokernel. So by 
the semipositivity of $\TT_{X'}|_{l'}$, we have 
$c_1(\LL)\cdot l=c_1(\LL')\cdot l'\geq 0$.
\end{proof}

As for the pseudo-effectivity of the second Chern classes, we may 
put forward the following

\begin{conj}\label{3.5}
Let $X$ be a variety of dimension $n$ with only terminal singularities 
and $-K_X$ nef. Then for any ample Cartier divisors $H_1,\cdots,H_{n-2}$ 
on $X$, we have 
\[ c_2(X)\cdot H_1\cdots H_{n-2}\geq 0. \]
\end{conj}

We define the numerical dimension $\nu(-K_X)$ of $-K_X$ to be the 
greatest nonnegative integer $\nu$ such that $(-K_X)^\nu\not\equiv 0$. 
If $\nu(-K_X)=0$, then Conjecture \ref{3.5} follows from the Miyaoka 
theorem (cf.\ \cite[Theorem 1.1]{mi}). 
If $\nu(-K_X)=1$, then Conjecture \ref{3.5} was proved in 
\cite[Corollary 6.2]{kmm04}. 
If $\nu(-K_X)=n$, then it is Theorem \ref{1.5}. 
While for $1<\nu(-K_X)<n$, Conjecture \ref{3.5} is still open.

For the three-dimensional case of Conjecture \ref{3.5}, we can use 
the Minimal Model Program to deal with the case when $\nu(-K_X)=2$ to 
obtain a partial result. We refer the reader to \cite{xie} for the 
details.

\end{document}